\documentclass[11pt,graphicx]{paper}

\usepackage{longtable,setspace}
\usepackage{graphicx,rotating,lscape} 
\usepackage{color}
\usepackage{amsthm}
\usepackage{float,epsfig,amsmath,amssymb,euscript,color}

\usepackage{authblk}

\usepackage{amsmath,amssymb,epsfig,capt-of,ifthen,calc}
\usepackage{float,euscript}

\usepackage{graphicx}
\usepackage{caption}
\usepackage{subcaption}

\usepackage{aeguill,aecompl}
\usepackage{color}

\usepackage{soul}
\usepackage{multirow}
\usepackage{dcolumn}
\usepackage{longtable}
\usepackage{supertabular}
\usepackage{setspace}

\newtheorem{conj}{Conjecture}

\definecolor{bleuclair}{rgb}{0.7, 0.7, 1.0}
\definecolor{rosepale}{rgb}{1.0, 0.7, 1.0}

\begin{document}

\title{Dealing with prime numbers I.: On the Goldbach conjecture}

\author{Fausto Martelli \thanks{faustom@princeton.edu}} 
\affil{Frick Chemistry Laboratory, Department of Chemistry, Princeton University, 08540 Princeton, New Jersey USA } 
\date{\today}

\maketitle

\abstract{In this paper we present some observations about the well-known \emph{Goldbach conjecture}. In particular we list and interpret some numerical results which allow us to formulate a relation between prime numbers and even integers. We can also determine very thin and low diverging ranges in which the probability of finding a prime is one.} 

\newpage

\section{Introduction}
In a letter dated $7$ June $1742$, the Prussian mathematician Christian Goldbach suggested to Leonhard Euler that every integer which 
can be written as the sum of two prime numbers, can also be written as the sum of as many primes as one wishes, 
until all terms are units. In the margin of the same letter he also proposed a second conjecture stating that every integer greater than $2$ can 
be written as the sum of three primes. 
Almost four weeks later, on $30$ June $1742$, Euler replied  in a letter and reminded of an earlier conversation they had, 
in which Goldbach remarked his original (and not marginal) conjecture: 
\begin{quote}Every even integer greater than 2 can be written as the sum of two primes.
\end{quote} 
In this letter  Euler pointed out:
\begin{quote}
Every even integer is a sum of two primes. 
I regard this as a completely certain theorem, although I cannot prove it.
\end{quote}
It is the worth recalling that Goldbach considered $1$ to be a prime number. Since this fact is not accepted any more nowadays, 
Goldbach's conjecture is formulated as 
\begin{conj}[SGC]\label{SGC}
Every even integer greater than or equal to $4$ can be written as a sum of two (odd) prime integers.
\end{conj}
This statement is called the  \emph{strong} Goldbach conjecture (SGC) in order  to distinguish it from  weaker corollaries. The SGC implies 
the conjecture that all odd numbers greater than $7$ are the sum of three odd primes, which is known today as the \emph{weak} or \emph{ternary} Goldbach 
conjecture. If the strong Goldbach conjecture is true, the weak Goldbach conjecture is true by implication.\\

Many progresses have been made in the last decade. In particular,
in 1997 has been shown that Goldbach conjecture is related to the generalized Riemann hypotesis in the sense that Riemann implies the Goldbach
weak conjecture for all numbers [1]. An extensive computational research has been also done on this direction [2].
In 2012 has been shown that every even number $n\ge4$ is in fact the sum of at most six primes, from which it follows that every odd number 
$n\ge5$ is the sum of at most seven primes, without using the Riemann Hypothesis [3] extending a previous result [4]. 
However, the biggest contribution has been recently obtained by Harald Helfgott who published a pair of papers claiming to improve 
major and minor arc estimates sufficiently to unconditionally prove the weak Goldbach conjecture [5,6]. \\

In this paper we take for granted the SGC and we conjecture something more. A formal proof of our conjecture would immediately
lead, by implication, to a proof of the Strong Goldbach Conjecture. We than present here our idea with numerical simulations which 
give an "experimental" proof up to $8\times10^{9}$.

\section{Considerations on the Goldbach conjecture}

According to SGC, we suppose that every even integer number $r\ge4$ can be written as the sum of two prime numbers $p$ and $p'$, not necessarily distinct:
\begin{equation}\label{SGC}\tag{SGC}
 r=p+p'.
\end{equation}

We introduce here our observations. Distribute the positive integers  into three infinite columns $C_1$, $C_2$ and $C_3$ in such a way that $C_1$ contains the numbers congruent to $1$ modulo $3$, $C_2$  the numbers congruent to $2$ modulo $3$ and  $C_3$ contains the multiples of $3$. Finally let $\alpha$ be the row index.

\begin{eqnarray}\label{eq:order1}
 \begin{array}{cccc}
  C_1 & C_2 & C_3 \\
  \\
  1 & 2 & 3 & \alpha=1\\
  4 & 5 & 6 & \alpha=2\\
  7 & 8 & 9 & \alpha=3\\
 10 & 11& 12& \alpha=4 \\
  \vdots & \vdots & \vdots 
  \end{array}
\end{eqnarray}

In what follows we will indicate every element of $\mathbb N^*$ expliciting its column and row indexes. More precisely $$t_n^\alpha\in\mathbb{N}^*$$ means the integer contained in the $n$-th $C_{n}$ column and in the $\alpha$-th row.  

In the same manner as above, we distribute the positive prime numbers in three columns in the following manner. Given the set of prime numbers $\mathbb{P}=\left\{2,3,5,7,\dots\right\}$, we order naturally them obtaining $p_1=2,\,p_2=3,\,p_3=5,\dots$ Then we put $p_k$ in the columns where $k$ was in (\ref{eq:order1}). We obtain:
\begin{eqnarray}\label{eq:order2}
 \begin{array}{cccc}
  C'_1 & C'_2 & C'_3 \\
  \\
  2 & 3 & 5 & \delta=1\\
  7 & 11& 13& \delta=2\\
  17& 19& 23& \delta=3\\
 29 & 31& 37& \delta=4 \\
  \vdots & \vdots & \vdots 
  \end{array}
\end{eqnarray}

As we did for integers, we will indicate every element of $\mathbb P$ expliciting its column and row indexes. More precisely $$p_m^\delta\in\mathbb{P}$$ means the integer contained in $C'_m$ and in the $\delta$-th row.  Observe that  every prime number has in this way two pairs of indexes: the former if it is seen as an element of $\mathbb{N}^*$, the latter if an element of $\mathbb{P}$. One must pay attention to not make confusion between them. When we work with integers  we use the symbol $t_n^\alpha$, otherwise we use $p_m^\delta$. 

\begin{figure}[h!]
   \centering
           \includegraphics[width=0.7\textwidth]{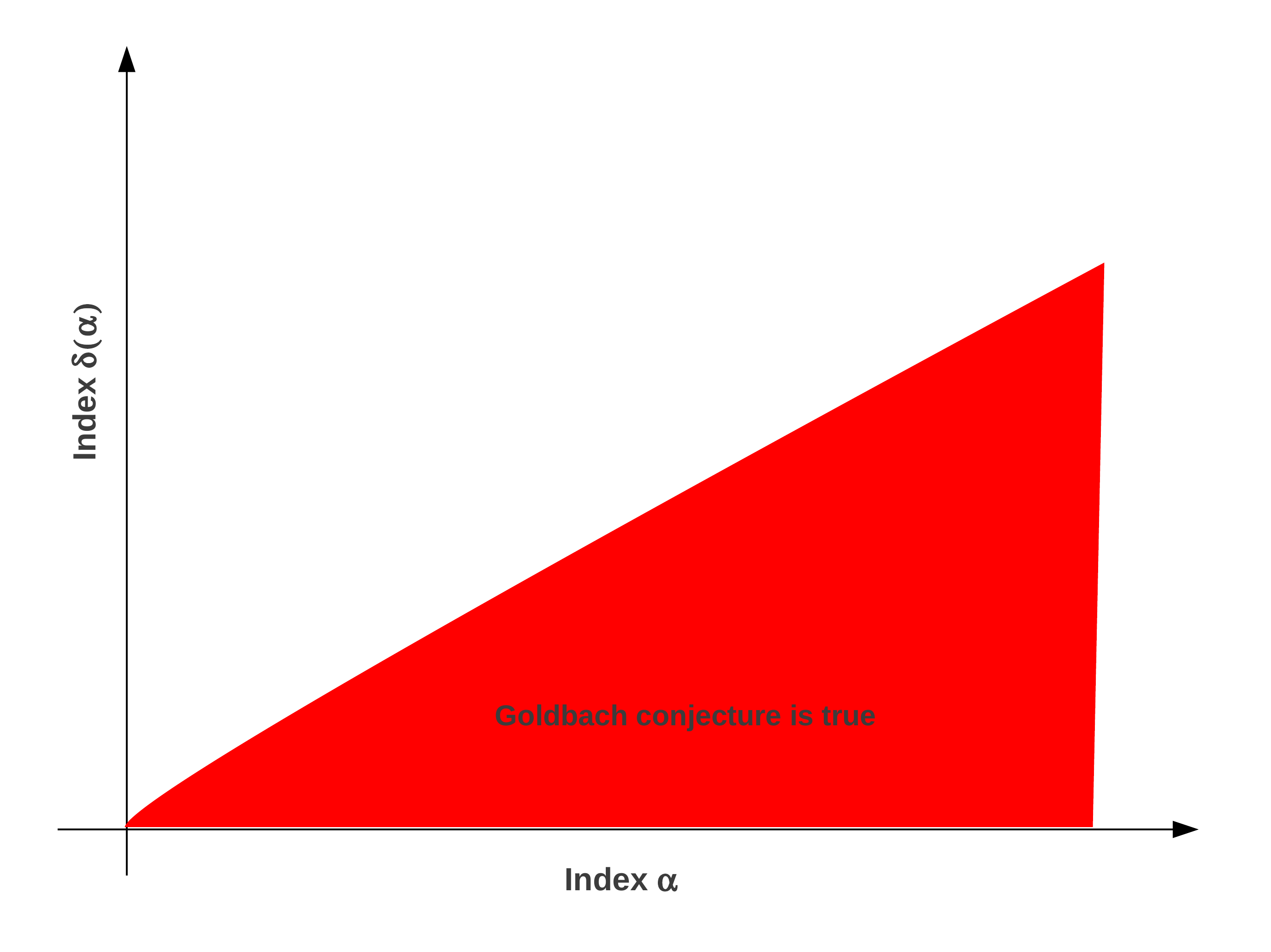}
           \caption{Space in which the Goldbach conjecture is valid (red).}
           \label{fig:spazi-1}
   \end{figure}

Now we are ready for the formulation of our conjecture:
\begin{conj}\label{sgc}
Let  $t_n^\alpha\in\mathbb{N}^*$ be a positive even integer. Then there esist two prime numbers $p_{n}^{\delta(\alpha)}, \,p_{m}^{\gamma}\in\mathbb{P}$ such that:
 \begin{equation}
  t_{n}^{\alpha}=p_{n}^{\delta(\alpha)}+p_{m}^{\gamma}  \label{eq:Gb2}
 \end{equation}
\end{conj}
where the the column index of one of the two primes is the same as the even number we decompose.\\  
After a fitting procedure on a large set, we have observed that the index $\delta(\alpha)$ follows Euler distribution:
\begin{equation}
 \delta(\alpha)=\left[B\frac{\alpha}{\log(\alpha)}\right] \label{eq:Euler}
\end{equation}
with $B$ constant. 
From a pictorial point of view, for large numbers the row indexes $\delta$ and $\alpha$ define a space in which Goldbace conjecture
is true (see fig(\ref{fig:spazi-1})). On the other hand, our conjecture allows us to select two thin slices
   \begin{figure}[h!]
   \centering
           \includegraphics[width=0.7\textwidth]{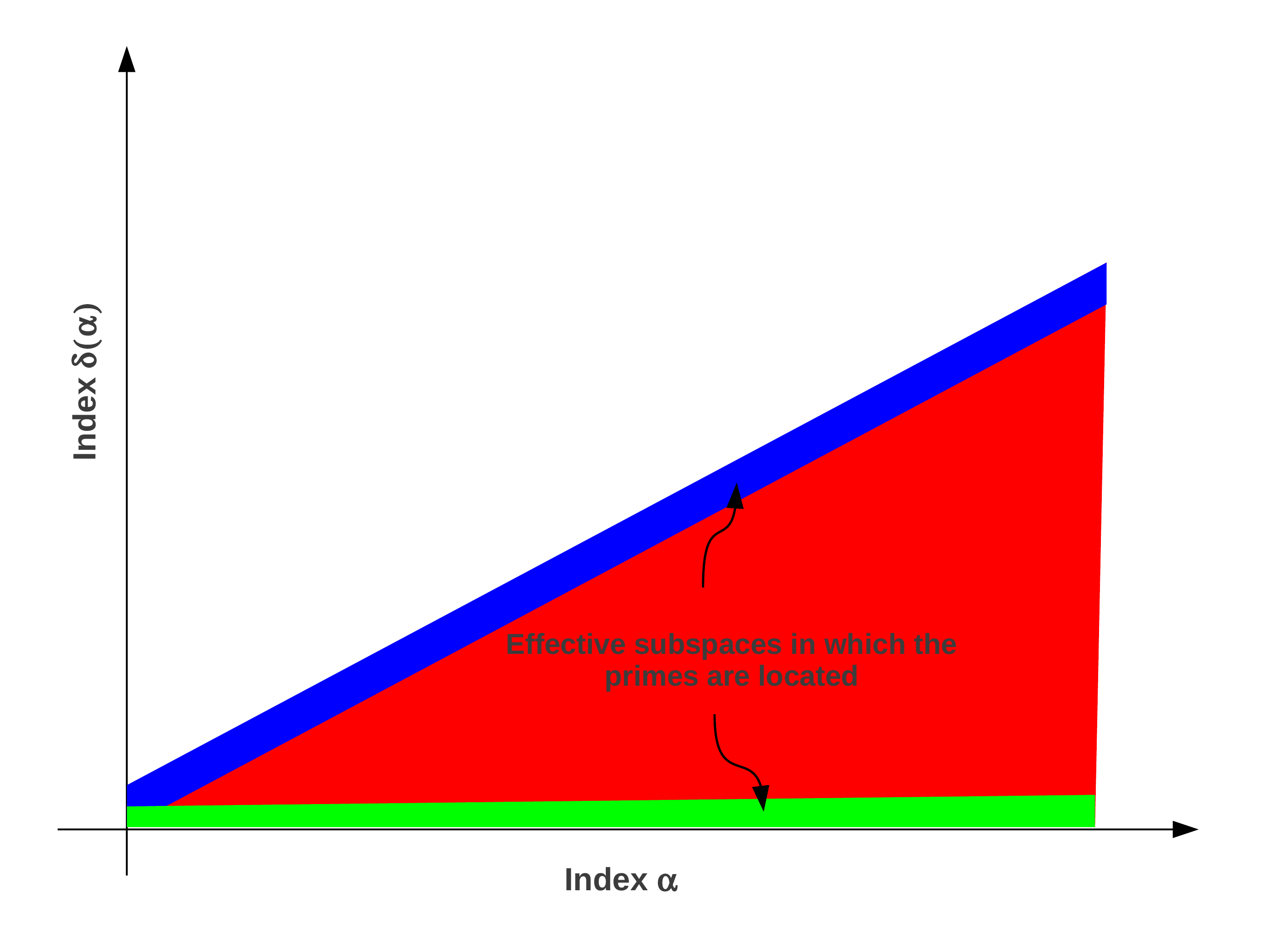}
           \caption{Slices in which our conjecture places the prime numbers satisfying Goldbach conjecture (blue and green).}
           \label{fig:spazi}
\end{figure}
of this space and to locate in there the prime numbers satisfying Goldbach (see fig.(\ref{fig:spazi})). In particular, $\delta$ define
the prime closer to the integer we want to examine, while $\gamma$ define the corresponding farthest, i.e., a small prime.

Although we do not have a formal proof of this result, we have performed numerical simulations in order to test our
conjecture which has been satisfiet for even integers up to 8433220000, whose closest prime is 8433219983 and the farthest is 17, 
apart for even integers 6, 16 and 164. 
Being these only three exceptions very small numbers, they must not be considered as a disproof of our conjecture. 
In fact, the distribution of $\delta(\alpha)$ in eq.(\ref{eq:Euler}) involving row indexes follow the prime number theorem
apart for a multiplicative factor. This point can be seen as a heuristic proof of our conjecture.

The multiplicative coefficient is $B=0.997602$. By multiplying Euler distribution, it minimize the distance 
\begin{equation}
 \left [ \delta(\alpha)\right ]-p_{n}^{\delta(\alpha)} 
\end{equation} 
with a correlation coefficient of 1.0000 between the calculated value for $\delta(\alpha)$ and the fitted value.

%

\section{Experimental results}
In this section we present our numerical results which give an "experimental" proof of our conjecture. We also briefly
present the program \textsf{CONJECTURE} and we explain the basic algorithm, but we do not go deep inside the computational performances
in order to do not distract the reader from the main message.
The program is freely available upon request to the author.\\
Numerical simulations has been performed on Princeton high-performance computers with
on a Red Hat 6, 1536 cores and a total RAM of 12 TB machine.
The multi-core platform processes theoretically on 16 Tflops with processor speed of 2.67 GHz Westmere.

\subsection{Numerical results}
\begin{figure}[h]
  \begin{center}
   \includegraphics[scale=.60]{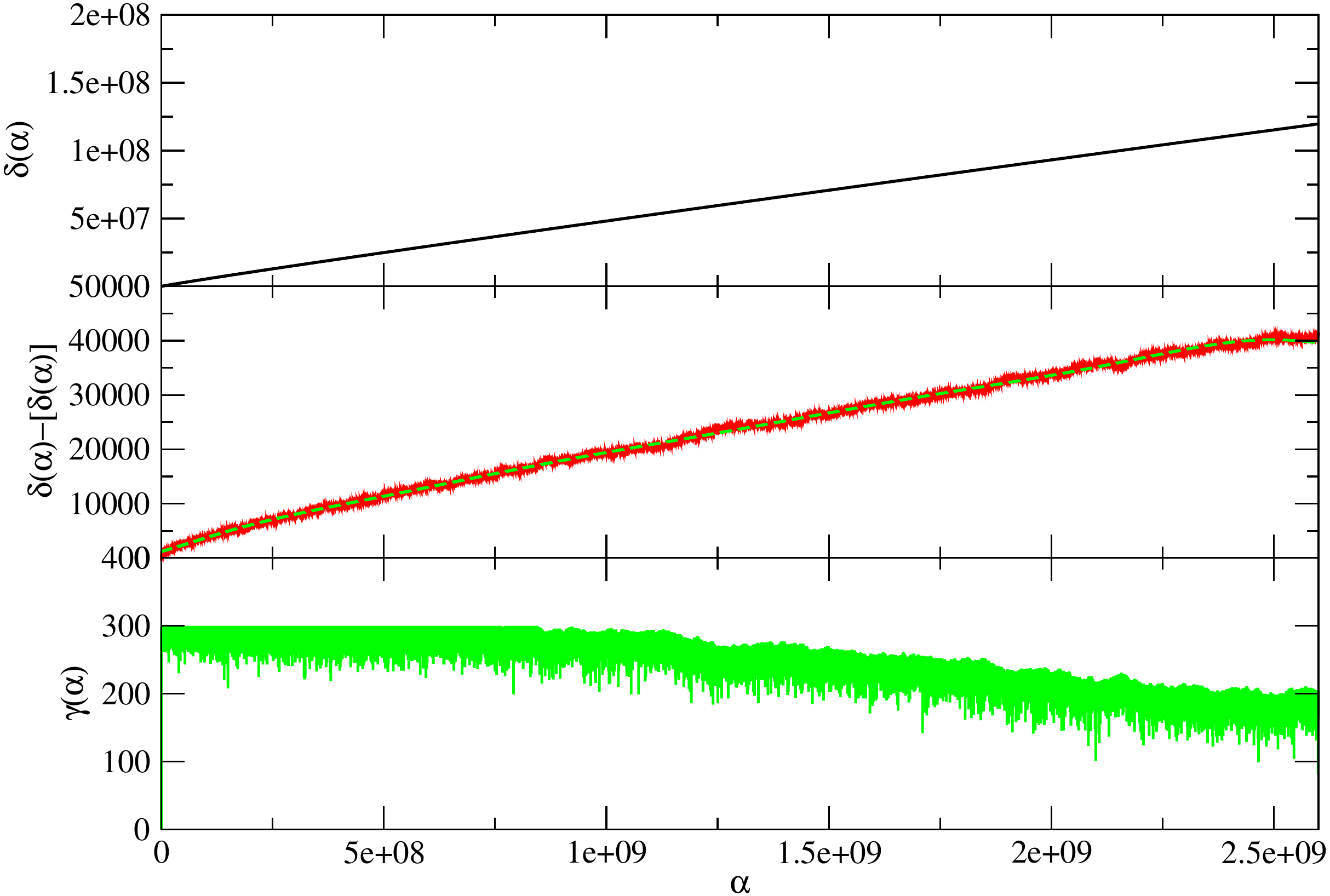}
   \caption{Upper panel: distribution of $\delta(\alpha)$ as a function of index $\alpha$. 
           Middle pane: continuous line represent the difference between 
           $\delta(\alpha)$ calculated by means of eq.(\ref{eq:Euler}) and its integer part; dashed line represent the fitted curve with
           regression procedure up to the 10-th order. Lower panel: distribution of $\gamma$ index.
           \label{fig:dist}}
 \end{center}
\end{figure}

\begin{figure}[h]
  \begin{center}
   \includegraphics[scale=.60]{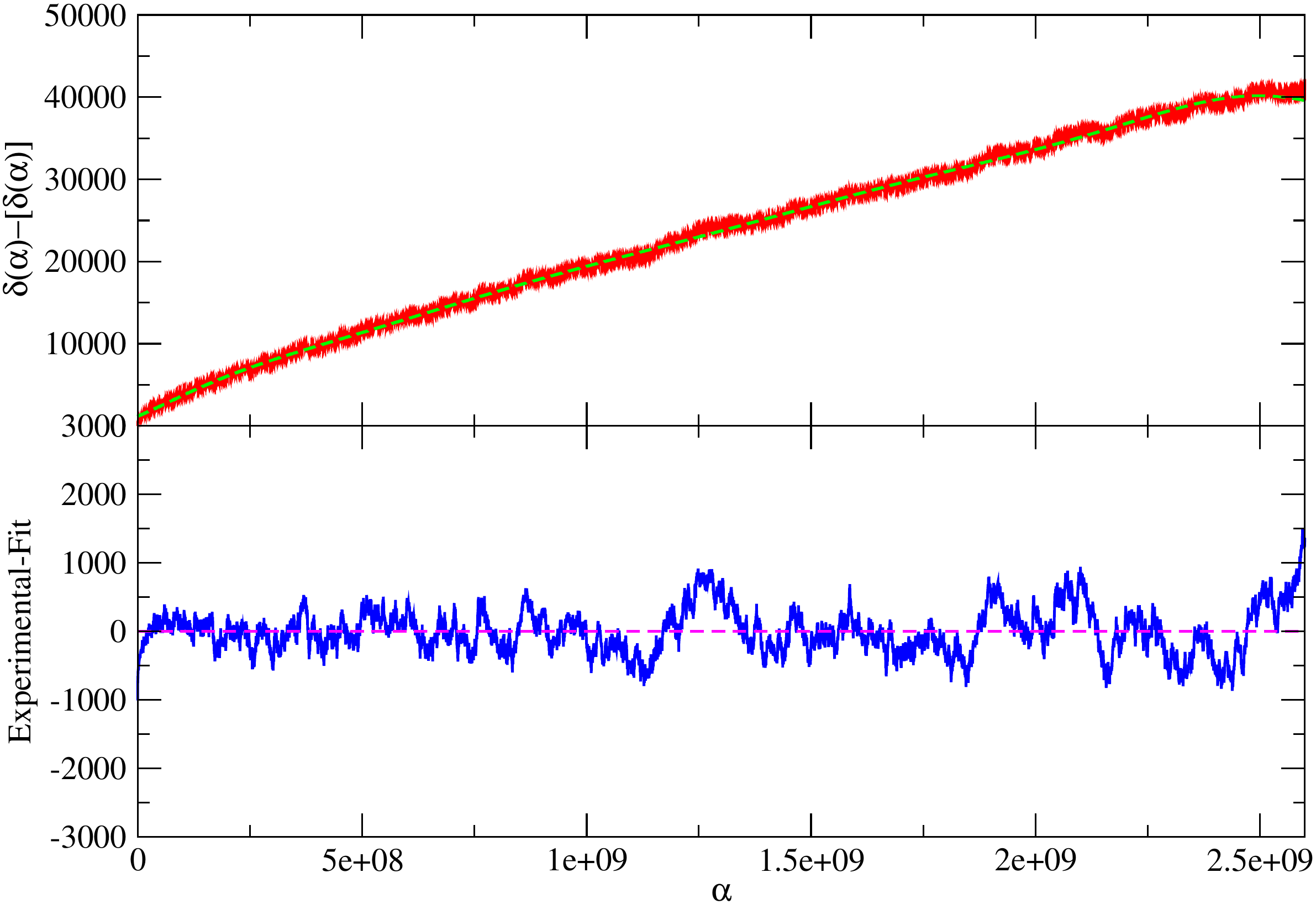}
   \caption{Upper panel: distribution of $\delta(\alpha)$ as a function of index $\alpha$. 
           Lower panel: difference between the experimental result and the fitted curve with the 10-th order polynomial 
           \label{fig:dist-2}}
 \end{center}
\end{figure}

We now present numerical results on our tests. In fig.(\ref{fig:dist}) we have three plots. 
The upper panel shows the distribusion $\delta(\alpha)$ as a function of $\alpha$. A fitting procedure of this curve with the 
Euler function shown in eq.(\ref{eq:Euler}) gives a perfect agreement with a correlation coefficient of 1.000000. 
In the middle panel the continuous curve shows the difference between $\delta(\alpha)$ calculated by means of 
eq.(\ref{eq:Euler}) and its integer part, and the dashed line represent the fitted curve with a regression up to the 10-th order.
The fitting polynomial is:
\begin{equation*}
 y=1184.9+2.251*10^{-5}x+5.68*10^{-14}x^{2}+...
\end{equation*}
where higher powers gives less significant contributions. The lower panel shows the distribution of $\gamma$. We can see that the
variation range is very thin. \\
The upper panel of figure (\ref{fig:dist-2}) shows again the the difference between $\delta(\alpha)$ calculated by means of
eq.(\ref{eq:Euler}), and its integer part.
The lower panel shows the difference between $\delta(\alpha)$ and the fitted values. We can see that these flactuations are
confined in a very thin range, showing than that the simulated curve very well reproduce the calculated values for $\delta(\alpha)$.

\subsection{The algorithm}
\begin{figure}[h]
  \begin{center}
   \includegraphics[scale=.60]{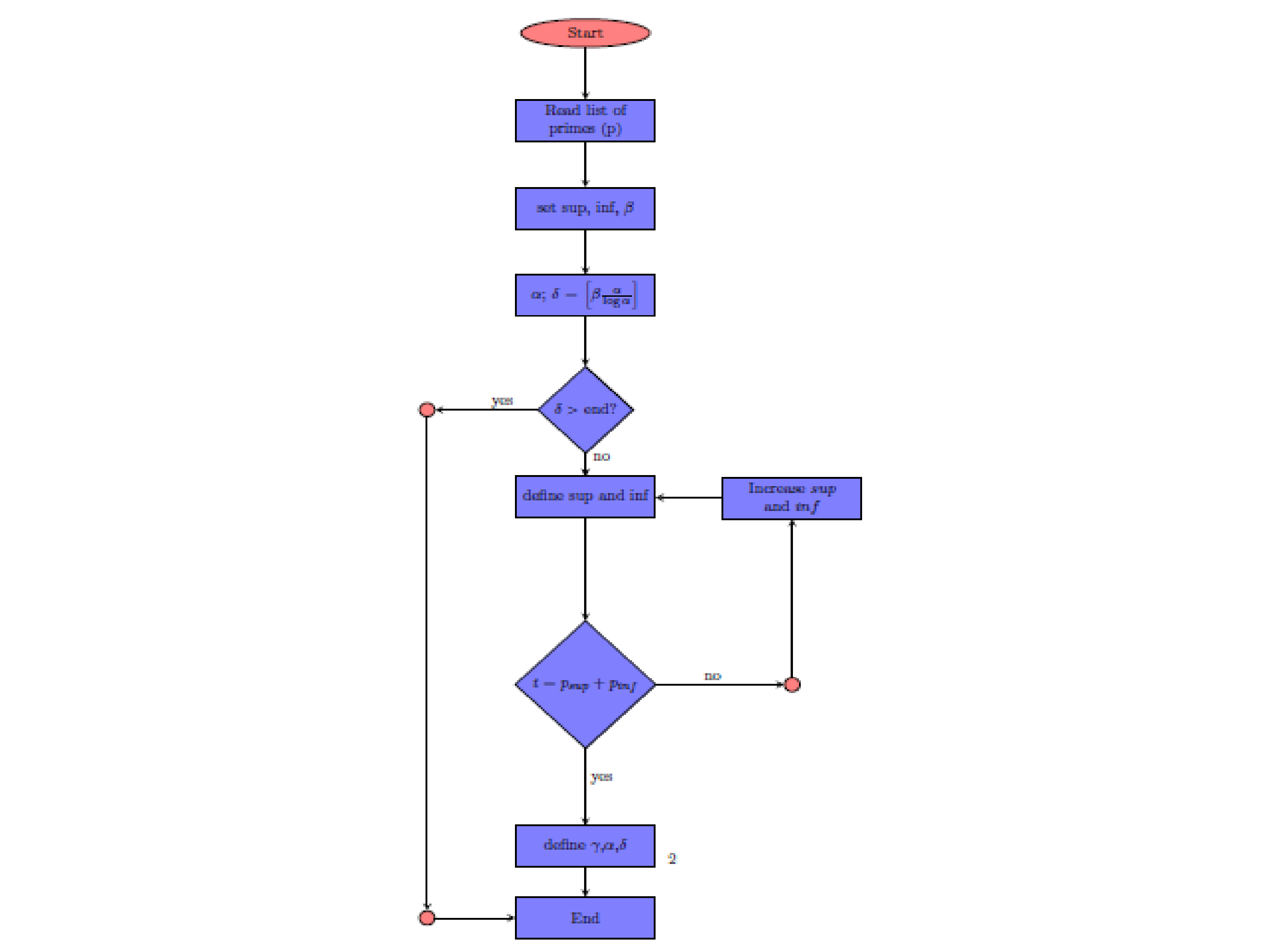}
   \caption{Basic flow chart of the \textsf{CONJECTURE} program.
           \label{fig:flow}}
 \end{center}
\end{figure}

The program \textsf{CONJECTURE} is a parallel program whose flow chart is reported in fig. (\ref{fig:flow}). 
In particular, of the first almost 9 billion integers we have split each billion on different nodes in order to
maximize the performances.\\
As input the program requires the
number $n$ representing the even natural up to which we want to test our conjecture, and a list of primes which
is stored in the file \textit{primes.dat}. The algorithm runs over $n$ and, at each step, evaluates when
$i=n (\mod 3)$, being $i$ the index over the naturals. When this relation is satisfied the program determines
$\alpha$, $\delta(\alpha)$ and $\left [\delta(\alpha)\right ]$, i.e., the index of the natural as defined
in eq.(\ref{eq:order1}), the real $\delta(\alpha)$ function and its integer part as defined in eq.(\ref{eq:Euler}) and 
representing the numerical function. \\
In order to minimize the computing time we define two ranges in which the algorithm looks for the primes. These two
ranges are self-adjustable, dependent of two parameters determined by fitting procedure:
\begin{equation*}
 \text{upper range:} [\delta-sup,\delta]\\
\end{equation*}
\begin{equation*}
 \text{lower range:} [\beta-inf, \beta+inf]
\end{equation*}
being $\delta$ defined in eq.(\ref{eq:Euler}), $\beta=$ and $sup$ and $inf$ initially posed as 10 and 80 respectively.
This two ranges correspond to the thin subspaces show pictorially on figs.(\ref{fig:spazi-1}),(\ref{fig:spazi}).\\
If the our conjecture is not satisfied, then the ranges are increased by increasing $sup$ and $inf$.
This step-adaptive characteristic increment enormously the performances of the program compared with the non 
step-adaptive version.

\section{Concluding remarks}
In this paper we have presented a conjecture whose formal proof would give a proof of the SGC. Numerical simulation performed up
to $8\times10^{9}$ do not show any exception apart for 6, 16 and 164. \\
Changing the number of columns on the natural and the prime number spaces (eqs.(\ref{eq:order1}),(\ref{eq:order2})) 
do not improve the results, suggesting
that the chosen number of column is the ideal way.

\section{References}
1. Deshouillers, Effinger, Te Riele and Zinoviev, "A complete Vinogradov 3-primes theorem under the Riemann hypothesis". Electronic Research Announcements of the American Mathematical Society, 1997 \textbf{3}, 99-104. \\
2. Saouter, "Checking the odd Goldbach Conjecture up to 10$^{20}$". Mathematics of Computation, 1998 \textbf{67}, 863-866.\\
3. Tao, "Every odd number greater than 1 is the sum of at most five primes". arXiv:1201.6656v4 (2012) \\
4. Kaniecki, "On Snirelman's constant under the Riemann hypothesis". Acta Arithmetica, 1995, \textbf{75}, 361-374.\\
5. Helfgott, "Major arcs for Goldbach's theorem". arXiv:1305.2897 (2013) \\
6. Helfgott, "Minor arcs for Goldbach's problem". arXiv:1205.5252 (2012) \\

\section*{Acknowledgements}
The author thanks Prof. Peter Sarnak (Institute for Advanced Studies of Princeton), Prof. Roberto Car (Princeton University),
Prof. Marie-Francoise Politis (Universit\'e 
Pierre et Marie Cure), Prof. Daniel Borgis and Rodolphe Vuilleumier (Ecole normale superieure of Paris) and
Dr. Riccardo Spezia (Universit\'e d'Evry val d'Essonne) for many useful discussions.

\end{document}